\documentclass[11pt]{amsart}

\usepackage{geometry}
\usepackage{colortbl}
\usepackage{amsmath,amssymb}
\usepackage{mathrsfs}
\usepackage{graphicx}
\usepackage{amsmath}
\usepackage{dsfont}
\usepackage{multirow}
\usepackage{graphicx}
\usepackage{epstopdf}
\usepackage{color}
\usepackage{ulem}
\usepackage{hyperref}
\usepackage{enumerate}
\usepackage{caption}
\usepackage{subfig}
\usepackage{bm}
\usepackage{overpic}

\usepackage{tikz}

\usepackage{varioref}
\usepackage{hyperref}
\usepackage{cleveref}

\makeatletter

\newcommand{\Rmnum}[1]{\expandafter\@slowromancap\romannumeral #1@}
\makeatother

\def\abs#1{\left|#1\right|}

\crefformat{figure}{Fig.~#2#1#3}
\crefmultiformat{figure}{Figures~(#2#1#3)}{ and~(#2#1#3)}{, (#2#1#3)}{ and~(#2#1#3)}
\crefformat{table}{Table~#2#1#3}

\DeclareMathOperator{\sech}{sech}

\theoremstyle{definition}

\theoremstyle{remark}

\newtheorem{remark}{Remark}

\numberwithin{equation}{section}
\begin{document}

\title{Enhanced second-order Gauss-Seidel projection methods for the Landau-Lifshitz equation}

\author{Panchi Li}
\address{Materials Innovation Institute for Life Sciences and Energy (MILES), HKU-SIRI, Shenzhen, China.}
\address{Department of Mathematics, The University of Hong Kong, Hong Kong, China}
\email{lipch@hku.hk}

\author{XIAO-PING WANG$^{\ast}$}\thanks{*Corresponding author}
\address{School of Science and Engineering, The Chinese University of Hong Kong, Shenzhen, Guangdong 518172, China}
\address{Shenzhen International Center for Industrial and Applied Mathematics, Shenzhen Research Institute of Big Data, Guangdong 518172, China}
\email{wangxiaoping@cuhk.edu.cn}

\graphicspath{{figures/}}


\date{\today}

\dedicatory{}
\keywords{Gauss-Seidel projection methods, Landau-Lifshitz equation, Backward differentiation formula, Second-order accuracy}
\begin{abstract}
The dynamics of magnetization in ferromagnetic materials are modeled by the Landau-Lifshitz equation, which presents significant challenges due to its inherent nonlinearity and non-convex constraint. These complexities necessitate efficient numerical methods for micromagnetics simulations. The Gauss-Seidel Projection Method (GSPM), first introduced in 2001, is among the most efficient techniques currently available. However, existing GSPMs are limited to first-order accuracy. This paper introduces two novel second-order accurate GSPMs based on a combination of the biharmonic equation and the second-order backward differentiation formula, achieving computational complexity comparable to that of solving the scalar biharmonic equation implicitly. The first proposed method achieves unconditional stability through Gauss-Seidel updates, while the second method exhibits conditional stability with a Courant-Friedrichs-Lewy constant of 0.25. Through consistency analysis and numerical experiments, we demonstrate the efficacy and reliability of these methods. Notably, the first method displays unconditional stability in micromagnetics simulations, even when the stray field is updated only once per time step.
\end{abstract}
\maketitle

\section{Introduction}
Ferromagnetic materials are ideal candidates for magnetic recording applications due to their intrinsic magnetic properties. This class of materials locally exhibits a net magnetic moment, known as magnetization, in the absence of external magnetic fields. The distribution of magnetization in the magnetic body gives rise to various static magnetic textures, such as the domain wall, magnetic vortex and skyrmions. When the magnetic field, spin-polarized current or temperature gradient is applied, the magnetic textures would be driven, collapsed and generated, where the dynamics of the magnetization is guided by the Landau-Lifshitz~(LL) equation~\cite{Gilbert1955,LandauLifshitz1935}.

The LL equation consists of the gyromagnetic and damping terms. The gyromagnetic term models the precession of the magnetization and conserves the magnetic free energy of the system, whereas the damping term governs the dissipation of energy. Below Curie's temperature, the length of magnetization is a constant at a point-wise level, and this infers the non-convex constraint in the LL equation. In addition, both the gyromagnetic term and damping term are strongly nonlinear. These accumulated features pose challenges in designing efficient and stable numerical methods for simulating the LL equation in micromagnetics simulations.

Over the past several decades, various numerical approaches have been developed for the micromagnetics simulations, see e.g. review papers~\cite{Abert2019,Cimrak2007,SIAMRev2006DevelopmentLLG} and references therein. The explicit methods, such as Runge-Kutta methods~\cite{10.1145/79505.79507,ROMEO2008464}, are favored in the early days. Explicit methods commonly suffer from stability where a small temporal step size is allowed, which results in the low efficiency of micromagnetics simulations even when certain adaptive strategies are adopted. On the contrary, the implicit Crank-Nicolson \cite{doi:10.1137/050631070,Fuwa2012} method maintains the length of magnetization and the energy law simultaneously. Nevertheless, the solvability of the implicit methods is strictly constrained with a ratio between the mesh size and time step size, and the iteration solver is required in applications due to the nonlinearity. These issues can be addressed by using linearized implicit methods, also known as the semi-implicit methods, such as the semi-implicit Euler method~\cite{10.1093/imanum/dri011}, implicit-explicit approaches~\cite{ALOUGES20121345}, the semi-implicit backward differentiation formula (BDF) methods~\cite{akrivis2021higher,xie2020second} and the semi-implicit Crank-Nicolson method~\cite{doi:10.1137/20M1335431}. Semi-implicit methods are characterized by their unconditional stability, unique solvability, and high efficiency. When a normalized step is further applied, the length of magnetization is preserved, although the energy law cannot be strictly maintained. Conversely, in linearized implicit methods, the preservation of magnetization length is not guaranteed, which may lead to inaccurate solution in long-time simulations. In this regard, high-order accurate methods exhibit their advantages, and some recent progress has been made such as the predictor-corrector methods \cite{doi:10.1137/22M1501143} and IMEX-RK method \cite{gui2024thirdorder}. In practice, the high-order methods not only provide a more accurate numerical solution but also can release the constriction of time step size. An intriguing example is provided in \cite{li2023micromagnetics}, where the use of the second-order BDF (BDF2) scheme permits a larger time step for simulating the LL equation with Dzyaloshinskii-Moriya interaction.

In addition to the aforementioned methods, there exists a class of numerical approaches known as Gauss-Seidel projection methods (GSPMs). These methods have a computational complexity comparable to that of solving the scalar heat equation implicitly and are faster than a 4th-order Runge-Kutta method \cite{10.1007/978-3-540-31852-1_17}. The GSPM was first proposed in 2001 \cite{WANG2001357} and has since been extensively used to study magnetization dynamics under various external controls. To enhance stability and computational efficiency, these methods have been developed in recent years \cite{garcia2003improved, Panchi2022DCDS-B, LI2020109046}. However, all available GSPMs are currently only first-order accurate in time.

In this paper, we present two second-order GSPMs, in which biharmonic equations are first solved to construct the second-order approximation of the Laplacian operator. By a combination of the semi-implicit BDF2 scheme, the second-order approximation for the LL equation is obtained. The computational complexity is therefore on the order of solving the scalar biharmonic equation implicitly. With the application of Gauss-Seidel updates, we numerically demonstrate one is unconditionally stable, while another is conditionally stable with the Courant-Friedrichs-Lewy (CFL) constant of 0.25. In addition, based on the proposed Scheme A we provide an unconditionally stable method for micromagnetics simulations. Numerical experiments demonstrate that it maintains unconditionally stable with respect to the damping parameter and reaches optimal efficiency in micromagnetics simulations.

The rest of the paper is organized as follows. In Section 2, we describe the fundamental model problem. We introduce the GSPM and present the second-order GSPMs in Section 3. Numerical experiments of 1D, 2D and 3D examples are conducted in Section 4. An application of the proposed Scheme A for micromagnetics simulations is presented in Section 5. The conclusions are drawn in Section 6.

\section{LL equation}
Let $\mathbf{M} = (M_1, M_2, M_3)^T$ denote the magnetization, and the LL equation reads as
\begin{equation}
  \partial_t\mathbf{M} = -\gamma\mathbf{M} \times \mathbf{H} - \frac{\gamma\alpha}{M_s} \mathbf{M} \times (\mathbf{M} \times \mathbf{H}),
  \label{equ:LL-equation-physics}
\end{equation}
where $\gamma$ is the gyromagnetic ratio, $0 < \alpha \ll 1$ is the dimensionless damping parameter, $M_s = |\mathbf{M}|$ is the saturation magnetization, and $\mathbf{H}$ is the effective magnetic field derived from the variational calculation of the energy functional
\begin{equation}
  E[\mathbf{M}] = \int_{\Omega} \left[\frac{A}{M_s^2}\abs{\nabla\mathbf{M}}^2 + \Phi\left(\frac{\mathbf{M}}{M_s}\right) - \mu_0\mathbf{H}_{\mathrm{e}}\cdot\mathbf{M} \right]\mathrm{d}\mathbf{x} + \frac{\mu_0}{2}\int_{\mathds{R}^3}\abs{\nabla U}^2\mathrm{d}\mathbf{x}.
  \label{equ:magnetic-free-energy}
\end{equation}
Here $\mu_0$ denotes the permeability of vacuum, $A$ is the exchange constant, and $\mathbf{H}_{\mathrm{e}}$ is the external magnetic field. In this paper, we will consider the uniaxial anisotropy with the easy-axis $\mathbf{e}_1 = (1, 0, 0)^T$, i.e., $\Phi\left(\mathbf{M}/M_s\right) = K_u(M_2^2+M_3^2)/(M_s^2)$ with $K_u$ the anisotropy constant. The last term on the right-hand side of \eqref{equ:magnetic-free-energy} is the dipolar energy from the stray field induced by the magnetization inhomogeneity and discontinuity, in which
\[
U = \int_{\Omega}\nabla N(\mathbf{x} - \mathbf{y})\cdot\mathbf{M}(\mathbf{y})\mathrm{d}\mathbf{y},
\]
where $N(\mathbf{x}) = -1/(4\pi \abs{\mathbf{x}})$. Hence the effective field is
\begin{equation}
  \label{equ:H_effective}
  \mathbf{H}_{\mathrm{eff}} = \frac{2A}{M_s^2}\Delta\mathbf{M} - \frac{2K_u}{ M_s^2}(M_2\mathbf{e}_2 + M_3\mathbf{e}_3) + \mu_0\mathbf{H}_{\mathrm{e}} +\mu_0 \mathbf{H}_{\mathrm{s}},
\end{equation}
where $\mathbf{e}_2 = (0, 1, 0)^T$, $\mathbf{e}_3 = (0, 0, 1)^T$ and $\mathbf{H}_{\mathrm{s}} = -\nabla U$ defines the stray field.

Using the definitions $\mathbf{H}_{\mathrm{eff}} = \mu_0 M_s\mathbf{h}$, $\mathbf{H}_{\mathrm{e}} =M_s\mathbf{h}_{\mathrm{e}}$, $\mathbf{H}_{\mathrm{s}} = M_s\mathbf{h}_{\mathrm{s}}$, and $\mathbf{M} = M_s\mathbf{m}$, with the spatial rescaling $\mathbf{x}\rightarrow L\mathbf{x}$ ($L$ the diameter of $\Omega$), the dimensionless form of the LL energy functional is
\begin{equation}
\label{equ:LLenergydimensionless}
\tilde{E}[\mathbf{m}] = \frac{1}{2}\int_{\Omega'}\left[\epsilon\abs{\nabla\mathbf{m}}^2 + q(m_2^2 + m_3^2) - 2\mathbf{h}_{\mathrm{e}}\cdot\mathbf{m} - 2\mathbf{h}_{\mathrm{s}}\cdot\mathbf{m}\right]\mathrm{d}\mathbf{x},
\end{equation}
where $\epsilon = 2A/(\mu_0M_s^2L^2)$, $q = 2K_u/(\mu_0M_s^2)$ and $\tilde{E}[\mathbf{m}] = {E}[\mathbf{M}]/(\mu_0M_s^2)$. Furthermore, with the temporal rescaling $t\rightarrow t(M_s\mu_0\gamma)^{-1}$, the dimensionless LL equation reads as
\begin{equation}
\label{equ:dimensionlessLLG}
\partial_t\mathbf{m} = -\mathbf{m}\times\mathbf{h} - \alpha\mathbf{m}\times\left( \mathbf{m} \times \mathbf{h} \right),
\end{equation}
where
\begin{equation}
  \label{equ:dimensionlesseffectivefield}
  \mathbf{h} = \epsilon\Delta\mathbf{m} - q(m_2\mathbf{e}_2 + m_3\mathbf{e}_3) + \mathbf{h}_{\mathrm{e}} + \mathbf{h}_{\mathrm{s}}.
\end{equation}
Over a bounded domain, the homogeneous Neumann boundary condition is used as
\begin{equation}
  \label{equ:Neumann}
  \partial_{\boldsymbol{\nu}}\mathbf{m} = 0,
\end{equation}
where $\boldsymbol{\nu}$ represents the outward unit normal vector along boundary surface $\partial\Omega$.

\section{GSPMs for the LL equation}
For the sake of simplicity and clarity, we consider the LL equation
\begin{equation}
  \partial_t\mathbf{m} = -\mathbf{m}\times \Delta\mathbf{m} - \alpha\mathbf{m}\times\left(\mathbf{m}\times\Delta\mathbf{m}\right),
  \label{equ:LL-equation}
\end{equation}
in which only the exchange field $\Delta \mathbf{m}$ is considered. The corresponding energy functional is
\begin{equation}
  E_e[\mathbf{m}] = \frac 12\int_{\Omega}|\nabla\mathbf{m}|^2 \mathrm{d}\mathbf{x}.
  \label{equ:LL-energy}
\end{equation}
Given a nonequilibrium initialization $\mathbf{m}_{\mathrm{in}}(\mathbf{x}) = \mathbf{m}_0(\mathbf{x})$ satisfying $|\mathbf{m}_0(\mathbf{x})| = 1$ in a point-wise sense, the energy \eqref{equ:LL-energy} decays as time evolves, and the length of magnetization is preserved in a point-wise sense simultaneously.

The earliest version of the GSPM was proposed in \cite{WANG2001357}, and subsequent improvements were made in \cite{garcia2003improved,LI2020109046}. An attempt to develop the second-order GSPM was made in \cite{Panchi2022DCDS-B}, where the GSPM-BDF2 enhanced the efficiency of micromagnetic simulations but did not achieve second-order accuracy. To explain the accuracy of this class of approaches clearly, we consider the gyromagnetic equation
\begin{equation}
  \partial_t\mathbf{m} = -\mathbf{m} \times \Delta\mathbf{m}.
  \label{equ:gyromagnetic-equation}
\end{equation}
The first-order GSPM solves it in fractional steps as follows:
\begin{align}
  \frac{\mathbf{m}^* - \mathbf{m}^n}{\Delta t} &= \Delta\mathbf{m}^*, \label{equ:1st-GSMP-step1}\\
  \frac{\mathbf{m}^{n+1} - \mathbf{m}^n}{\Delta t} &= -\mathbf{m}^n \times \frac{\mathbf{m}^* - \mathbf{m}^n}{\Delta t}.\label{equ:1st-GSMP-step2}
\end{align}
It is evident that only heat equations need to be solved during the time stepping, which allows the Fast Fourier Transform technique to be conveniently applied when the regularity domain is considered.

The first-order approximation of \eqref{equ:1st-GSMP-step2} for the gyromagnetic equation \eqref{equ:gyromagnetic-equation} can be obtained by the Taylor expansion. According to the equation \eqref{equ:1st-GSMP-step1}, we have
\begin{equation}
  \mathbf{m}^* = (I - \Delta t\Delta)^{-1}\mathbf{m}^n = \mathbf{m}^n + \Delta t\Delta\mathbf{m}^n + \Delta t^2\Delta^2\mathbf{m}^n + \cdots.
  \label{equ:Taylor-m-star}
\end{equation}
A combination of \eqref{equ:1st-GSMP-step2} and \eqref{equ:Taylor-m-star} yields
\begin{equation}
  \frac{\mathbf{m}^{n+1} - \mathbf{m}^n}{\Delta t} = -\mathbf{m}^n \times \frac{\mathbf{m}^* - \mathbf{m}^n}{\Delta t} = -\mathbf{m}^n \times \Delta\mathbf{m}^n - \Delta t\mathbf{m}^n \times \Delta^2\mathbf{m}^n + \cdots,
  \label{equ:Taylor-gyromagnetic}
\end{equation}
which indicates that if $\mathbf{m} \in C^1([0, T]; [C^0(\overline\Omega)]^3) \cap L^{\infty}([0, T]; [C^4(\overline\Omega)]^3)$, the numerical scheme \eqref{equ:1st-GSMP-step1}-\eqref{equ:1st-GSMP-step2} shall provide the first-order accurate solution in time. Obviously, \eqref{equ:Taylor-gyromagnetic} is conditionally stable, and the stability can be improved by the application of the Gauss-Seidel update for components of magnetization.

Similar to the GSPM-BDF2 as in \cite{Panchi2022DCDS-B}, the construction of second-order GSPM is based on the semi-implicit BDF2 scheme as
\begin{equation}
  \frac{\frac 32\mathbf{m}^{n+1} - 2\mathbf{m}^n + \frac 12\mathbf{m}^{n-1}}{\Delta t} = -\hat{\mathbf{m}}^{n+1} \times \Delta \mathbf{m}^{n+1},
  \label{equ:semi-implicit-BDF2-gyro-equ}
\end{equation}
where $\hat{\mathbf{m}}^{n+1} = 2\mathbf{m}^n - \mathbf{m}^{n-1}$. According to \eqref{equ:Taylor-m-star} and \eqref{equ:Taylor-gyromagnetic}, we see that $\mathbf{m}^*$ is the first-order approximation of $\mathbf{m}^n$, as well as $\frac{\mathbf{m}^* - \mathbf{m}^n}{\Delta t}$ for $\Delta\mathbf{m}^n$. Hence the direct combination of \eqref{equ:Taylor-m-star} and \eqref{equ:semi-implicit-BDF2-gyro-equ}, GSPM-BDF2 method, only achieves the first-order accuracute approximation.

Instead of the Taylor expansion \eqref{equ:Taylor-m-star}, we consider
\begin{equation}
  (I - \Delta t\Delta + \Delta t^2\Delta^2)^{-1} = I + \Delta t\Delta - \Delta t^3\Delta^3 + \cdots.
\end{equation}
Therefore, we solve the biharmonic equation with the backward Euler method, and the second-order GSPM is
\begin{gather}
  \frac{\mathbf{m}^{*} - \hat{\mathbf{m}}^{n+1}}{\Delta t} = \Delta\mathbf{m}^{*} - \Delta t\Delta^2\mathbf{m}^{*}, \label{equ:2st-GSMP-step1}\\
  \frac{\frac 32\mathbf{m}^{n+1} - 2\mathbf{m}^n + \frac 12\mathbf{m}^{n-1}}{\Delta t} = -\hat{\mathbf{m}}^{n+1} \times \frac{\hat{\mathbf{m}}^* - \hat{\mathbf{m}}^{n+1}}{\Delta t}.\label{equ:2st-GSMP-step2}
\end{gather}
Using a similar expansion as in \eqref{equ:Taylor-gyromagnetic}, we have
\begin{align}
  \frac{\frac 32\mathbf{m}^{n+1} - 2\mathbf{m}^n + \frac 12\mathbf{m}^{n-1}}{\Delta t} &= -\hat{\mathbf{m}}^{n+1} \times \frac{\hat{\mathbf{m}}^{n+1} + \Delta t\Delta\hat{\mathbf{m}}^{n+1} - \Delta t^3\Delta^3\hat{\mathbf{m}}^{n+1} + \cdots - \hat{\mathbf{m}}^{n+1}}{\Delta t} \nonumber \\
  &= -\hat{\mathbf{m}}^{n+1} \times \Delta\hat{\mathbf{m}}^{n+1} - \Delta t^2\hat{\mathbf{m}}^{n+1} \times\Delta^3\hat{\mathbf{m}}^{n+1} + \cdots.
  \label{equ:Taylor-expansion-2}
\end{align}
This indicates that if $\mathbf{m} \in C^3([0, T]; [C^0(\overline\Omega)]^3) \cap C^2([0, T]; [C^6(\overline\Omega)]^3)$, the numerical scheme \eqref{equ:2st-GSMP-step1}-\eqref{equ:2st-GSMP-step2} shall provide the second-order accurate solution.

With the damping term being added, here we provide two GSPMs with second-order accuracy for the LL equation \eqref{equ:LL-equation}. Both of our methods are based on the following semi-implicit BDF2 discretization
\begin{equation}
  \frac{\frac 32\mathbf{m}^{n+1} - 2\mathbf{m}^n + \frac 12\mathbf{m}^{n-1}}{\Delta t} = -\hat{\mathbf{m}}^{n+1}\times\Delta\mathbf{m}^{n+1} - \alpha\hat{\mathbf{m}}^{n+1}\times(\hat{\mathbf{m}}^{n+1}\times \Delta\mathbf{m}^{n+1}).
  \label{equ:semi-implicit-BDF2}
\end{equation}
We substitute \eqref{equ:2st-GSMP-step1} into the above \eqref{equ:semi-implicit-BDF2}, and the full LL equation is approximated by
\begin{gather}
  \frac{\mathbf{m}^* - \hat{\mathbf{m}}^{n+1}}{\Delta t} = \Delta \mathbf{m}^* - \Delta t\Delta^2\mathbf{m}^{*}, \label{equ:SI-step1}\\
  \frac{\frac 32\tilde{\mathbf{m}}^{n+1} - 2\mathbf{m}^n + \frac 12\mathbf{m}^{n-1}}{\Delta t} = -\hat{\mathbf{m}}^{n+1} \times \frac{\mathbf{m}^* - \hat{\mathbf{m}}^{n+1}}{\Delta t} - \alpha\hat{\mathbf{m}}^{n+1} \times \left(\hat{\mathbf{m}}^{n+1} \times \frac{\mathbf{m}^* - \hat{\mathbf{m}}^{n+1}}{\Delta t}\right)\label{equ:SI-step2}
\end{gather}
with the projection $ \mathbf{m}^{n+1} = \tilde{\mathbf{m}}^{n+1} / |\tilde{\mathbf{m}}^{n+1}|$ in a point-wise sense.
In following, the formula
$$\hat{\mathbf{m}}^{n+1}\times(\hat{\mathbf{m}}^{n+1}\times \Delta\mathbf{m}^{n+1}) = (\hat{\mathbf{m}}^{n+1}\cdot \Delta\mathbf{m}^{n+1}) \hat{\mathbf{m}}^{n+1} - |\hat{\mathbf{m}}^{n+1}|^2 \cdot \Delta\mathbf{m}^{n+1}$$
will be used, in which $|\hat{\mathbf{m}}^{n+1}|$ is maintained because of $|\hat{\mathbf{m}}^{n+1}| \neq 1$.

\begin{remark}
  The continuous counterpart of the equation \eqref{equ:SI-step1} is
  \begin{equation}
    \partial_t \mathbf{m} = \Delta \mathbf{m} + \tau\Delta^2\mathbf{m}.
    \label{equ:continuous-biharmonic-equation}
  \end{equation}
  where $\tau$ is a small parameter to be determined. This defines an auxiliary equation to obtain the second-order approximation for the Laplacian operator. According to the above analysis, we set $\tau = -\Delta t$ in the semi-discrete scheme. To ensure the well-posedness of the auxiliary equation \eqref{equ:continuous-biharmonic-equation}, one more boundary condition related to the second-order spatial derivative $\partial_{\nu}^2\mathbf{m}$ is required. Here we consider $\partial_{\nu}^2\mathbf{m} = 0$, and then more ghost points around boundaries would be introduced to discrete the bihormonic operator.
\end{remark}

Using the Taylor expansion for the BDF2 scheme \eqref{equ:semi-implicit-BDF2}, we can see that $\mathbf{m}^{n+1} = \mathbf{m}(t_{n+1}) + \mathcal{O}(\Delta t^2)$, where $t_{n+1} = (n+1)\Delta t$. Similarly, we also have $\tilde{\mathbf{m}}^{n+1} = \mathbf{m}(t_{n+1}) + \mathcal{O}(\Delta t^2)$ in \eqref{equ:SI-step2}. Take $\tilde{\mathbf{m}}^{n+1} = \mathbf{m}(t_{n+1}) + \mathbf{M}\Delta t^2$ with $\mathbf{M}$ being the coefficient of the second-order term, and we then obtain
\begin{align*}
  |\tilde{\mathbf{m}}^{n+1}|^2 = |\mathbf{m}(t_{n+1})|^2 + 2\mathbf{m}(t_{n+1}) \cdot \mathbf{M}\Delta t^2 + |\mathbf{M}|^2\Delta t^4
  \leq (1+\epsilon) + \left(1+\frac 1{\epsilon}\right)|\mathbf{M}|^2\Delta t^4,
\end{align*}
and
\begin{equation*}
  |\tilde{\mathbf{m}}^{n+1}|^2 \geq (1-\epsilon) + \left(1 - \frac 1{\epsilon}\right)|\mathbf{M}|^2\Delta t^4.
\end{equation*}
Assume that there exist constants $C_0$ and $\epsilon_0 \in (-\epsilon, \epsilon)$ such that $|\tilde{\mathbf{m}}^{n+1}| = (1 + \epsilon_0) + C_0\Delta t^2$. We arrive at
\begin{align*}
  \mathbf{m}^{n+1} = \frac{\tilde{\mathbf{m}}^{n+1}}{|\tilde{\mathbf{m}}^{n+1}|} = \frac{\mathbf{m}(t_{n+1}) + \mathcal{O}(\Delta t^2)}{1 + \epsilon_0 + C\Delta t^2} = (1 - \epsilon_0)\mathbf{m}(t_{n+1}) + \mathcal{O}(\Delta t^2),
\end{align*}
where $\epsilon_0$ is an arbitrarily small constant. To ensure $|\mathbf{m}^{n+1}| = 1$, it is obvious that we can set $\epsilon_0$ to be the order of $\Delta t^2$. The second-order approximation is therefore maintained. An example of the rigorous error estimate for the projection of second-order accuracy can be found in \cite{doi:10.1137/20M1335431}.

\subsection{Scheme A}
Denote $\mathcal{L} = I - \Delta t\Delta + \Delta t^2\Delta^2$. With the Gauss-Seidel update for the components of magnetization, we solve the full LL equation in two steps:
\begin{enumerate}[Step 1:]
  \item semi-implicit Gauss-Seidel step
  \begin{equation*}
    \begin{aligned}
      g_i^{*} &= \mathcal{L}^{-1}\hat{m}_i^{n+1}, i = 1,2,3,\qquad   g_i^{n+1} = \mathcal{L}^{-1}\hat{m}_i^{*}, i = 1,2, \\
      \frac 32m_1^* &= 2{m}_1^n - \frac 12{m}_1^{n-1} -(\hat{m}_2^{n+1}g_3^* - \hat{m}_3^{n+1}g_2^*)
      \\&\qquad - \alpha (\hat{m}_1^{n+1}g_1^* + \hat{m}_2^{n+1}g_2^* + \hat{m}_3^{n+1}g_3^*)\hat{m}_1^{n+1} + \alpha ((\hat{m}_1^{n+1})^2 + (\hat{m}_2^{n+1})^2 + (\hat{m}_2^{n+1})^2) g_1^*, \\
      \frac 32m_2^* &= 2{m}_2^n - \frac 12{m}_2^{n-1} -(\hat{m}_3^{n+1}g_1^{n+1} - \hat{m}_1^{*}g_3^*)
      \\&\qquad - \alpha (\hat{m}_1^{*}g_1^{n+1} + \hat{m}_2^{n+1}g_2^* + \hat{m}_3^{n+1}g_3^*)\hat{m}_2^{n+1} + \alpha ((\hat{m}_1^{*})^2 + (\hat{m}_2^{n+1})^2 + (\hat{m}_2^{n+1})^2) g_2^*, \\
      \frac 32m_3^* &= 2{m}_3^n - \frac 12{m}_3^{n-1} -(\hat{m}_1^{*}g_2^{n+1} - \hat{m}_2^{*}g_1^{n+1})
      \\&\qquad - \alpha (\hat{m}_1^*g_1^{n+1} + \hat{m}_2^*g_2^{n+1} + \hat{m}_3^{n+1}g_3^*)\hat{m}_3^{n+1} + \alpha ((\hat{m}_1^*)^2 + (\hat{m}_2^*)^2 + (\hat{m}_2^{n+1})^2) g_3^*,
    \end{aligned}
  \end{equation*}
  where $\hat{m}_i^{*} = 2m_i^* - m_i^n (i = 1, 2)$.
  \item Projection onto $\mathcal{S}^2$
  \begin{equation*}
    \begin{pmatrix}
      m_1^{n+1} \\ m_2^{n+1} \\ m_3^{n+1}
    \end{pmatrix} = \frac 1{|\mathbf{m}^*|}\begin{pmatrix}
      m_1^{*} \\ m_2^{*} \\ m_3^{*}
    \end{pmatrix}.
  \end{equation*}
\end{enumerate}

In this scheme, the biharmonic equation with constant coefficients is solved 5 times at each time step. It is the second-order counterpart of Scheme A in \cite{LI2020109046}.

\subsection{Scheme B}
Following the similar approach of Scheme B in \cite{LI2020109046}, here we present the corresponding second-order version. This requires additional initial conditions $g_i^0 ( i = 1, 2, 3 )$. In practice, we first compute the numerical solution $\mathbf{m}_h^1$ using the first-order GSPM. The initial conditions are then determined by
\begin{equation*}
  g_i^0 = \mathcal{L}^{-1}\left( 2m_i^1 - m_i^0 \right), \quad i = 1, 2, 3.
\end{equation*}
The formal algorithm is as follows.
\begin{enumerate}[Step 1:]
  \item Semi-implicit Gauss-Seidel step
  \begin{align*}
    g_i^{n+1} &= \mathcal{L}^{-1}\hat{m}_i^{*}, i = 1,2,3, \\
    \frac 32m_1^* &= 2{m}_1^n - \frac 12{m}_1^{n-1} -(\hat{m}_2^{n+1}g_3^{n} - \hat{m}_3^{n+1}g_2^{n})
      \\&\qquad - \alpha (\hat{m}_1^{n+1}g_1^{n} + \hat{m}_2^{n+1}g_2^{n} + \hat{m}_3^{n+1}g_3^{n})\hat{m}_1^{n+1} + \alpha ((\hat{m}_1^{n+1})^2 + (\hat{m}_2^{n+1})^2 + (\hat{m}_2^{n+1})^2) g_1^{n} \\
      \frac 32m_2^* &= 2{m}_2^n - \frac 12{m}_2^{n-1} -(\hat{m}_3^{n+1}g_1^{n+1} - \hat{m}_1^{*}g_3^{n})
      \\&\qquad - \alpha (\hat{m}_1^{*}g_1^{n+1} + \hat{m}_2^{n+1}g_2^{n} + \hat{m}_3^{n+1}g_3^{n})\hat{m}_2^{n+1} + \alpha ((\hat{m}_1^{*})^2 + (\hat{m}_2^{n+1})^2 + (\hat{m}_2^{n+1})^2) g_2^{n} \\
      \frac 32m_3^* &= 2{m}_3^n - \frac 12{m}_3^{n-1} -(\hat{m}_1^{*}g_2^{n+1} - \hat{m}_2^{*}g_1^{n+1})
      \\&\qquad - \alpha (\hat{m}_1^*g_1^{n+1} + \hat{m}_2^*g_2^{n+1} + \hat{m}_3^{n+1}g_3^{n})\hat{m}_3^{n+1} + \alpha ((\hat{m}_1^*)^2 + (\hat{m}_2^*)^2 + (\hat{m}_2^{n+1})^2) g_3^{n}
  \end{align*}
  \item Projection onto $\mathcal{S}^2$
  \begin{equation*}
    \begin{pmatrix}
      m_1^{n+1} \\ m_2^{n+1} \\ m_3^{n+1}
    \end{pmatrix} = \frac 1{|\mathbf{m}^*|}\begin{pmatrix}
      m_1^{*} \\ m_2^{*} \\ m_3^{*}
    \end{pmatrix}.
  \end{equation*}
\end{enumerate}
In this scheme, the biharmonic equation with constant coefficients is solved 3 times at each time step.

\subsection{Spatial discretization}
Without loss of generality, we consider the uniform mesh for $\Omega \subset \mathds{R}^3$ with the mesh size $h = \Delta x = \Delta y = \Delta z$. Define the index $i = -1, 0, 1, \cdots, N_x + 1, N_x + 2$, $j = -1, 0, 1, \cdots, N_y+1, N_y+2$ and $k = -1, 0, 1, \cdots, N_z+1, N_z+2$, where $i = -1, 0, N_x+1, N_x+2$, $j = -1, 0, N_y+1, N_y+2$ and $k = -1, 0, N_z+1, N_z+2$ denote the "ghost points". We use the half grid points with $\mathbf{m}_{i,j,k} = \mathbf{m}((i-\frac 12)h, (j-\frac 12)h, (k-\frac 12)h)$. Consequently, the discrete Laplacian operator is
\begin{align*}
  \Delta_h \mathbf{m}_{i,j,k} = &\frac{\mathbf{m}_{i+1,j,k} -2\mathbf{m}_{i,j,k} + \mathbf{m}_{i-1,j,k}}{\Delta x^2} + \\
  &\frac{\mathbf{m}_{i,j+1,k} -2\mathbf{m}_{i,j,k} + \mathbf{m}_{i,j-1,k}}{\Delta y^2} + \\
  &\frac{\mathbf{m}_{i,j,k+1} -2\mathbf{m}_{i,j,k} + \mathbf{m}_{i,j,k-1}}{\Delta z^2}
\end{align*}
with $i = 1, \cdots, N_x$, $j = 1, \cdots, N_y$ and $k = 1, \cdots, N_z$. The values of "ghost points" are determined by the homogeneous boundary condition such that
\begin{gather*}
  \mathbf{m}_{0, j, k} = \mathbf{m}_{1, j, k}, \quad \mathbf{m}_{N_x+1, j, k} = \mathbf{m}_{N_x, j, k},\\
  \mathbf{m}_{i, 0, k} = \mathbf{m}_{i, 1, k}, \quad \mathbf{m}_{i, N_y+1, k} = \mathbf{m}_{i, N_y, k},\\
  \mathbf{m}_{i, j, 0} = \mathbf{m}_{i, j, 1}, \quad \mathbf{m}_{i, j, N_z+1} = \mathbf{m}_{i, j, N_z}.
\end{gather*}
In addition, since the biharmonic operator is introduced, the discrete biharmonic operator is
\begin{align*}
  \Delta_h^2 \mathbf{m}_{i,j,k} =
  &\frac{\Delta_h\mathbf{m}_{i+1,j,k} -2\Delta_h\mathbf{m}_{i,j,k} + \Delta_h\mathbf{m}_{i-1,j,k}}{\Delta x^2} + \\
  &\frac{\Delta_h\mathbf{m}_{i,j+1,k} -2\Delta_h\mathbf{m}_{i,j,k} + \Delta_h\mathbf{m}_{i,j-1,k}}{\Delta y^2} + \\
  &\frac{\Delta_h\mathbf{m}_{i,j,k+1} -2\Delta_h\mathbf{m}_{i,j,k} + \Delta_h\mathbf{m}_{i,j,k-1}}{\Delta z^2}.
\end{align*}
On the boundaries, we further introduce
\begin{gather*}
  \Delta_h \mathbf{m}_{0, j, k} = \Delta_h \mathbf{m}_{1, j, k}, \quad \Delta_h \mathbf{m}_{N_x+1, j, k} = \Delta_h \mathbf{m}_{N_x, j, k},\\
  \Delta_h \mathbf{m}_{i, 0, k} = \Delta_h \mathbf{m}_{i, 1, k}, \quad \Delta_h \mathbf{m}_{i, N_y+1, k} = \Delta_h \mathbf{m}_{i, N_y, k},\\
  \Delta_h \mathbf{m}_{i, j, 0} = \Delta_h \mathbf{m}_{i, j, 1}, \quad \Delta_h \mathbf{m}_{i, j, N_z+1} = \Delta_h \mathbf{m}_{i, j, N_z}.
\end{gather*}
These indicate that the remaining "ghost points" satisfy
\begin{gather*}
  \mathbf{m}_{-1, j, k} = \mathbf{m}_{2, j, k}, \quad \mathbf{m}_{N_x+2, j, k} = \mathbf{m}_{N_x-1, j, k},\\
  \mathbf{m}_{i, -1, k} = \mathbf{m}_{i, 2, k}, \quad \mathbf{m}_{i, N_y+2, k} = \mathbf{m}_{i, N_y-1, k},\\
  \mathbf{m}_{i, j, -1} = \mathbf{m}_{i, j, 2}, \quad \mathbf{m}_{i, j, N_z+2} = \mathbf{m}_{i, j, N_z-1},
\end{gather*}
which can be obtained by the homogeneous boundary condition as well. Applying the Taylor expansion, it is easy to check the second-order approximation of the above discrete forms for both Laplacian and biharmonic operators.

\section{Numerical experiments}
In this section, we conduct numerical experiments to validate the second-order approximation of the proposed GSPMs. The LL equation with a source term induced by an exact solution is considered in both 1D and 3D tests, while the LL equation without the source term is considered in a 2D test. We find that Scheme B is conditionally stable, and hence its CFL constant is investigated numerically.

\subsection{Numerical accuracy}
\subsubsection{1D case} We assume that there is a source $\mathbf{g}$ satisfying the LL equation
\begin{equation}
  \partial_t\mathbf{m} = -\mathbf{m} \times \Delta\mathbf{m} - \alpha\mathbf{m} \times (\mathbf{m} \times \Delta\mathbf{m}) + \mathbf{g}.
  \label{equ:LL-euq-with-source}
\end{equation}
The source $\mathbf{g}$ is the result of the exact solution
\begin{equation}
  \mathbf{m}(x, t) = (\cos(\bar{x})\sin(t), \sin(\bar{x})\sin(t), \cos(t))^T,
  \label{equ:1d-exact-solution}
\end{equation}
where $\bar{x} = x^2(1-x)^2$. We set $\Omega = (0, 1)$ and $\alpha = 0.01$ and check the numerical convergence rates in both time and space. The convergence rate with respect to time step size $\Delta t$ is shown in \Cref{tab:1D-convergence-rate-time}.
\begin{table}[htbp]
  \centering
  \caption{The numerical convergence rate of the two GSPMs in time. We fix $\alpha = 0.01$, $\Delta x = $1.0e-04 and $T = 0.3$.}
  \begin{tabular}{||c|cccc|c||}
    \hline
    $\Delta t$ & T/200 & T/300 & T/400 & T/500 & order \\
    \hline
    Scheme A & 1.5771e-04 & 7.4962e-05 & 3.9885e-05 & 2.3881e-05 & 2.06 \\
    Scheme B & 1.5754e-04 & 7.4324e-05 & 3.9032e-05 & 2.4842e-05 & 2.03 \\
    \hline
  \end{tabular}
  \label{tab:1D-convergence-rate-time}
\end{table}

Meanwhile, we set the time step size to $\Delta t = $1.0e-05, and vary the spatial mesh size with $\Delta x = 0.2, 0.1, 0.05, 0.04$. The second-order convergence rate in space is shown in \Cref{tab:1d-convergence-rate-in-space}.
\begin{table}[htbp]
  \centering
  \caption{The numerical convergence rate of the two GSPMs in space. We fix $\alpha = 0.01$, $\Delta t = $1.0e-05 and $T = 0.05$.}
  \begin{tabular}{||c|cccc|c||}
    \hline
    $\Delta x$ & 0.2 & 0.1 & 0.05 & 0.04 & order \\
    \hline
    Scheme A & 1.7327e-03 & 4.4315e-04 & 1.1888e-04 & 8.0147e-05 & 1.91 \\
    Scheme B & 1.7309e-03 & 4.4277e-04 & 1.1763e-04 & 7.8625e-05 & 1.92 \\
    \hline
  \end{tabular}
  \label{tab:1d-convergence-rate-in-space}
\end{table}

\subsubsection{2D case without source}
\begin{figure}[htbp]
  \centering
  \includegraphics[width = 3.in]{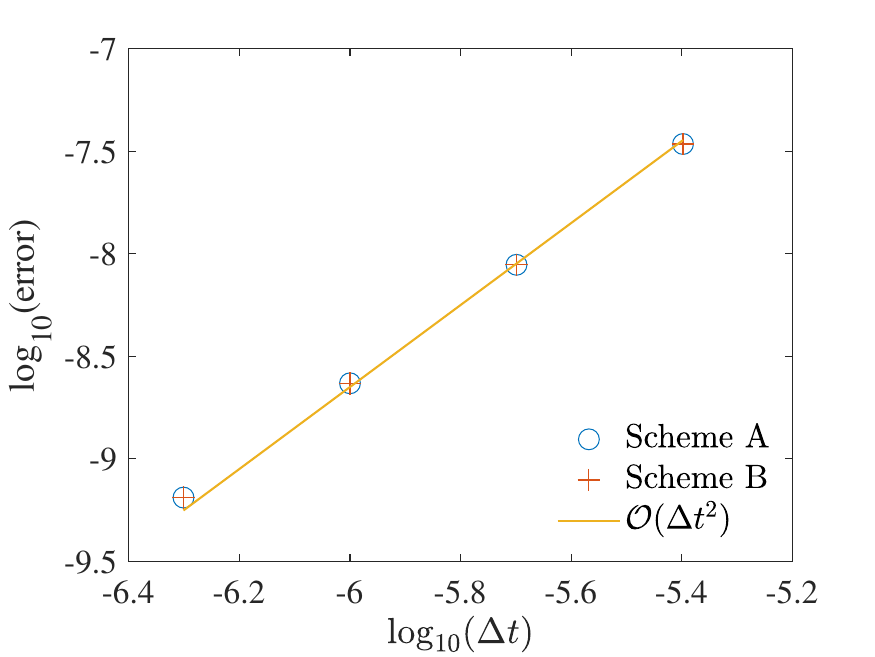}
  \caption{The convergence rate of Scheme A and Scheme B for the full LL equation without source term. The slope of the fitted solid line is 1.91.}
  \label{fig:2d-numerical-convergence}
\end{figure}

In this experiment, over the domain $\Omega = (0, 1)\times(0, 0.2)$, we adopt the initial condition
\begin{equation}
  \mathbf{m}_0(\mathbf{x}) = (\tanh(\ell), \sech(\ell), 0)^T,
  \label{equ:Neel-wall-initilization}
\end{equation}
where $\ell = (0.5-x)/ 2\eta$ with $\eta$ being a small quality. Here we set $\eta = \Delta x$ and $\alpha = 0.01$, and employ the semi-implicit BDF2 projection method \cite{xie2020second} to calculate the reference solution with $\Delta t = T/5000$. Other setups are: $\Delta x = \Delta y = 0.1$, $T = $4.0e-05. We vary the time step size $\Delta t = T/$[10 20 40 80], and obtain the second-order convergence rate as shown in \Cref{fig:2d-numerical-convergence}.

\subsubsection{3D case}
Consider the 3D equation with the source term induced by the exact solution
\begin{equation*}
  \mathbf{m}(\mathbf{x}, t) = (\cos(\bar x\bar y\bar z)\sin(t), \sin(\bar x\bar y\bar z)\sin(t), \cos(t))^T,
\end{equation*}
where $\bar y = y^2(1-y)^2$ and $\bar z = z^2(1-z)^2$. We fix the damping parameter $\alpha = 0.1$, the terminal time $T = 4.0$, and the computational domain $\Omega = (0, 1)^3$ with uniform mesh size $h = \Delta x = \Delta y = \Delta z$. The second-order convergence rates in both space and time are obtained as in~\Cref{tab:3d-convergence-space} and \Cref{tab:3d-convergence-time}.

\begin{table}[htbp]
  \centering
  \caption{The numerical convergence rate in space ($\Delta t = $1.0e-03).}
  \begin{tabular}{||c|cccc|c||}
    \hline
    $h$ & 1/6 & 1/8 & 1/10 & 1/12 & order \\
    \hline
    Scheme A & 1.1510e-04 & 6.4244e-05 & 4.0743e-05 & 2.7983e-05 & 2.04 \\
    Scheme B & 1.1619e-04 & 6.4855e-05 & 4.1130e-05 & 2.8272e-05 & 2.04 \\
    \hline
  \end{tabular}
  \label{tab:3d-convergence-space}
\end{table}

\begin{table}[htbp]
  \centering
  \caption{The numerical convergence rate in time ($h = 0.1$).}
  \begin{tabular}{||c|cccc|c||}
    \hline
    $\Delta t$ & 2.0e-01 & 1.0e-01 & 8.0e-02 & 5.0e-02 & order \\
    \hline
    Scheme A & 3.5853e-02 & 9.5834e-03 & 6.1886e-03 & 2.4294e-03 & 1.94 \\
    Scheme B & 3.3793e-02 & 9.3001e-03 & 6.0396e-03 & 2.3909e-03 & 1.91 \\
    \hline
  \end{tabular}
  \label{tab:3d-convergence-time}
\end{table}

\subsection{The stability of second-order GSPMs}
We investigate the stability constraint of the proposed methods in this experiment. The 1D LL equation \eqref{equ:LL-euq-with-source} with the exact solution \eqref{equ:1d-exact-solution} is employed. We fix $\alpha = 1.0$ and $T = 1.0$, and check the requirement of time step size for both Scheme A and Scheme B as in \Cref{tab:CFL-condition}. Numerical results show that Scheme B is constrained by the standard CFL condition with a constant of 0.25, while Scheme A is unconditional.
\begin{table}[htbp]
  \centering
  \caption{The stability constraint of two schemes.}
  \begin{tabular}{||c|c|c||}
    \hline
    & Scheme A & Scheme B \\
    \hline
    $h = 0.1$ & & $\Delta t \leq $2.5e-03 \\
    $h = 0.01$ & No restriction & $\Delta t \leq $2.5e-05\\
    $h = 0.005$ & & $\Delta t \leq $6.25e-06\\
    $h = 0.002$ & & $\Delta t \leq $1.0e-06 \\
    \hline
  \end{tabular}
  \label{tab:CFL-condition}
\end{table}

\section{An unconditionally stable second-order GSPM for micromagnetics simulations}
We further adopt Scheme A to check the stability with respect to damping parameter $\alpha$ in micromagnetics simulations. In following, when the full LL equation \eqref{equ:dimensionlessLLG} is considered, we split the effective field \eqref{equ:dimensionlesseffectivefield} into
\begin{equation}
  \mathbf{h} = \epsilon\Delta\mathbf{m} + \mathbf{f}(\mathbf{m})
\end{equation}
with $\mathbf{f}(\mathbf{m}) = - q(m_2\mathbf{e}_2 + m_3\mathbf{e}_3) + \mathbf{h}_{\mathrm{e}} + \mathbf{h}_{\mathrm{s}}$.
We therefore solve the equation
\begin{equation}
  \partial_t\mathbf{m} = -\mathbf{m}\times( \epsilon\Delta\mathbf{m} + \mathbf{f}(\mathbf{m}) ) - \alpha\mathbf{m}\times\left( \mathbf{m} \times (\epsilon\Delta\mathbf{m} + \mathbf{f}(\mathbf{m})) \right)
\end{equation}
in two steps as follows.
\begin{enumerate}[Step 1:]
  \item semi-implicit Gauss-Seidel step
  \begin{gather*}
    g_i^{*} = (I - \epsilon\Delta t\Delta_h + \epsilon^2\Delta t^2\Delta_h^2)^{-1}(\hat{m}_i^{n+1} + \Delta tf_i(\hat{m}_i^{n+1})), i = 1,2,3,\\
    g_i^{n+1} = (I - \epsilon\Delta t\Delta_h + \epsilon^2\Delta t^2\Delta_h^2)^{-1}(\hat{m}_i^{*} + \Delta tf_i(\hat{m}_i^{n+1})), i = 1,2,
  \end{gather*}
  \begin{equation*}
    \begin{aligned}
      \frac 32m_1^* &= 2{m}_1^n - \frac 12{m}_1^{n-1} -(\hat{m}_2^{n+1}g_3^* - \hat{m}_3^{n+1}g_2^*)
      \\&\qquad - \alpha (\hat{m}_1^{n+1}g_1^* + \hat{m}_2^{n+1}g_2^* + \hat{m}_3^{n+1}g_3^*)\hat{m}_1^{n+1} + \alpha ((\hat{m}_1^{n+1})^2 + (\hat{m}_2^{n+1})^2 + (\hat{m}_2^{n+1})^2) g_1^*, \\
      \frac 32m_2^* &= 2{m}_2^n - \frac 12{m}_2^{n-1} -(\hat{m}_3^{n+1}g_1^{n+1} - \hat{m}_1^{*}g_3^*)
      \\&\qquad - \alpha (\hat{m}_1^{*}g_1^{n+1} + \hat{m}_2^{n+1}g_2^* + \hat{m}_3^{n+1}g_3^*)\hat{m}_2^{n+1} + \alpha ((\hat{m}_1^{*})^2 + (\hat{m}_2^{n+1})^2 + (\hat{m}_2^{n+1})^2) g_2^*, \\
      \frac 32m_3^* &= 2{m}_3^n - \frac 12{m}_3^{n-1} -(\hat{m}_1^{*}g_2^{n+1} - \hat{m}_2^{*}g_1^{n+1})
      \\&\qquad - \alpha (\hat{m}_1^*g_1^{n+1} + \hat{m}_2^*g_2^{n+1} + \hat{m}_3^{n+1}g_3^*)\hat{m}_3^{n+1} + \alpha ((\hat{m}_1^*)^2 + (\hat{m}_2^*)^2 + (\hat{m}_2^{n+1})^2) g_3^*.
    \end{aligned}
  \end{equation*}
  \item Projection onto $\mathcal{S}^2$
  \begin{equation*}
    \begin{pmatrix}
      m_1^{n+1} \\ m_2^{n+1} \\ m_3^{n+1}
    \end{pmatrix} = \frac 1{|\mathbf{m}^*|}\begin{pmatrix}
      m_1^{*} \\ m_2^{*} \\ m_3^{*}
    \end{pmatrix}.
  \end{equation*}
\end{enumerate}
Notably, in the above algorithm, the magnetic field $\mathbf{f}(\mathbf{m})$ is not updated in the Gauss-Seidel manner, i.e., the stray field is only calculated once at each time step. Therefore, we primarily assess the stability of the method with respect to the small damping parameter. We consider a thin-film ferromagnet with the size of $1\times 1\times 0.02 \;\mu\mathrm{m}^3$ and a cell size $0.004\;\mu\mathrm{m}^3$. The physical constants characteristic of the ferromagnet are: $A = 1.3\times 10^{-11}\;\mathrm{J/m}$, $M_s = 8 \times 10^{5} \;\mathrm{A/m}$, $K_u = 1\times 10^2 \;\mathrm{J/m^3}$ and $\gamma = 1.76 \times 10^{11} \mathrm{T}^{-1}\mathrm{s}^{-1}$. In simulations, we use the initial state $\mathbf{m}_0 = (0, 1, 0)^T$ if $x \in [0, L_x/5] \cup [4L_x/5, L_x]$ and $\mathbf{m}_0 = (1, 0, 0)^T$ otherwise.

\begin{figure}[htbp]
  \centering
  \subfloat[Angle profile, magnetization profile and energy decaying of the full LL equation with $\alpha = 0.1$.]{\includegraphics[width=2.2in]{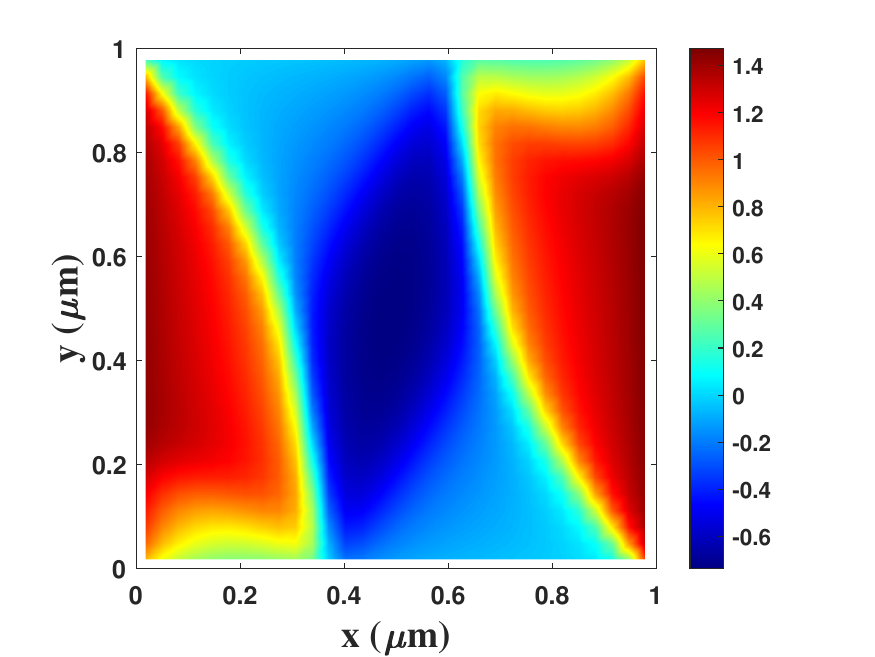}
  \includegraphics[width=2.2in]{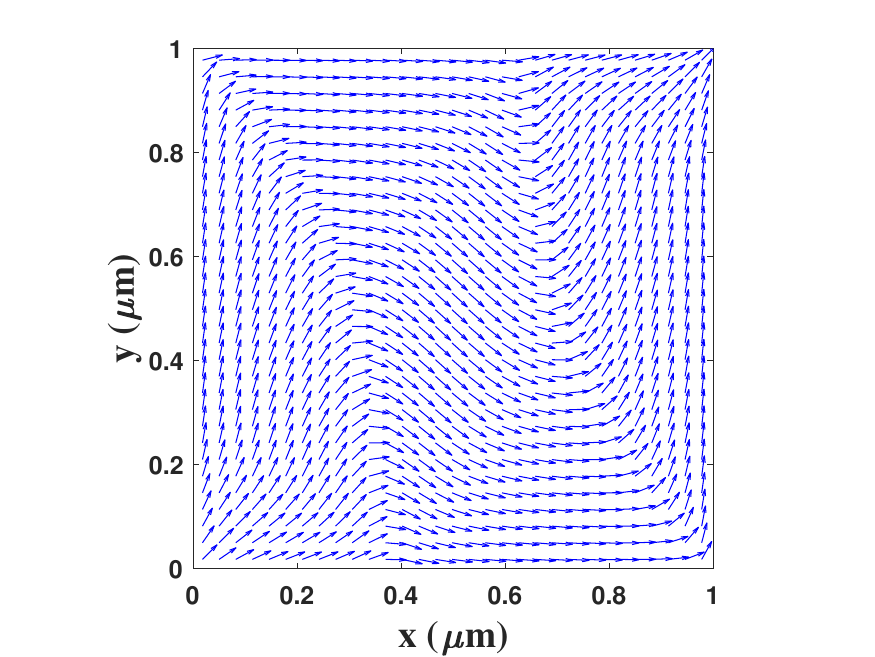}
  \includegraphics[width=2.2in]{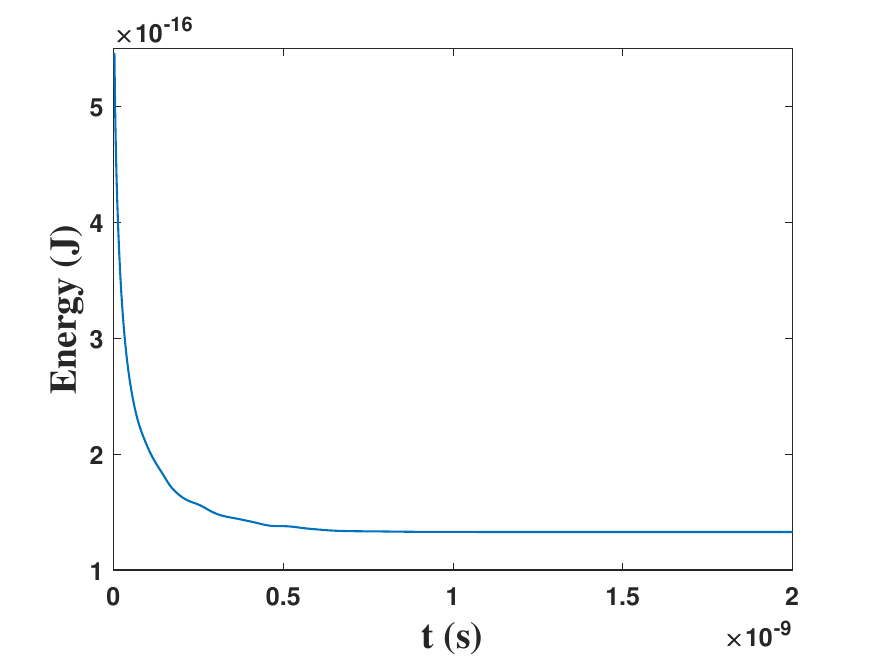}}\\
  \subfloat[Angle profile, magnetization profile and energy decaying of the full LL equation with $\alpha = 0.01$.]{\includegraphics[width=2.2in]{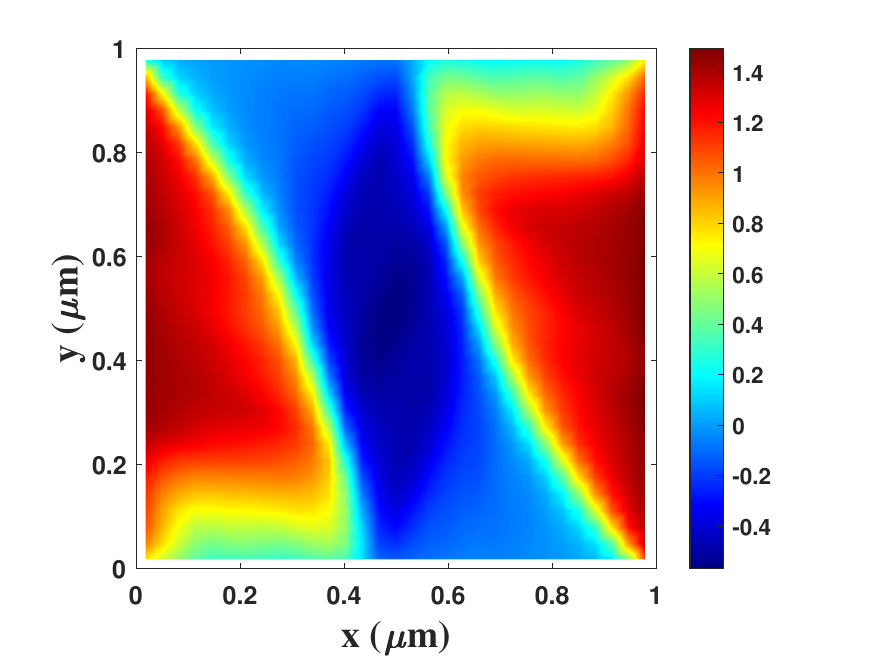}
  \includegraphics[width=2.2in]{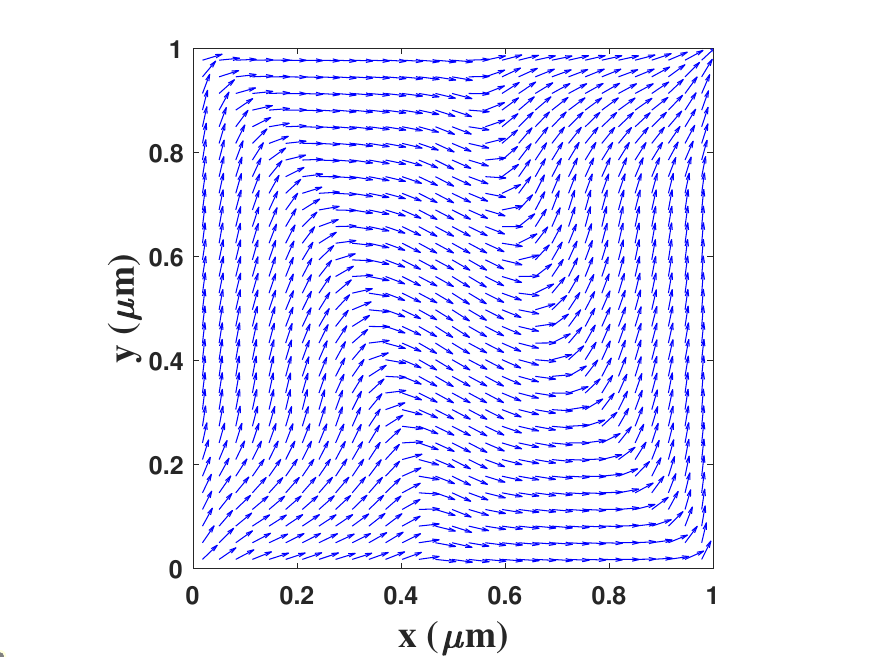}
  \includegraphics[width=2.2in]{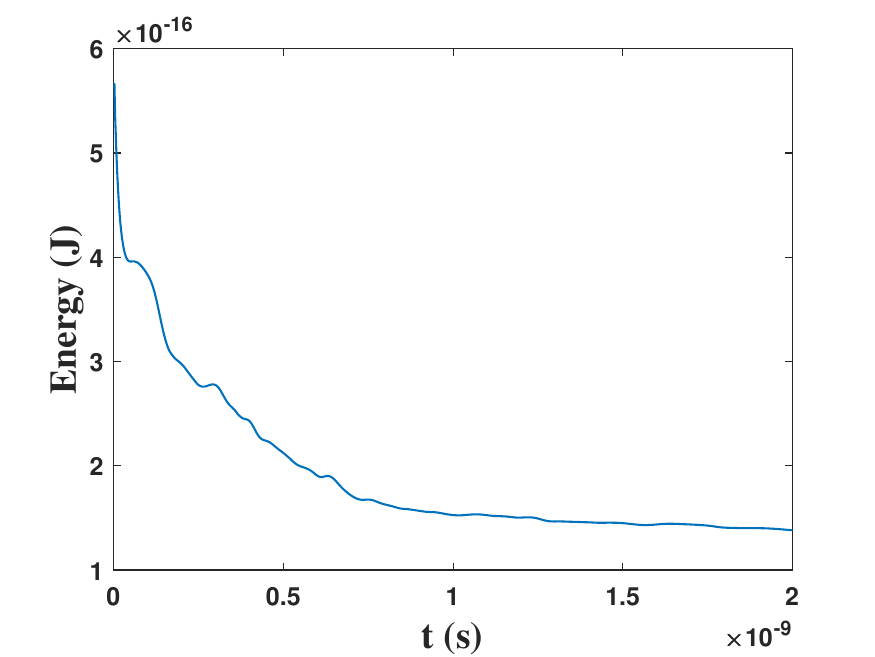}}
  \caption{The simulation of the full LL equation without external fields. Here the magnetization on the centered slice of the material in the $xy$ plane is depicted.}
  \label{fig:macromagnetics-simulations}
\end{figure}

The GSPMs provided in \cite{garcia2003improved,LI2020109046,WANG2001357} suffer from the computations of the stray field. Inadequate updates for the stray field would lead to the stability issue. To improve the computational efficiency and maintain stability, BDF2 was adopted in the GSPM-BDF2 where the stray field updates only once per time step, and then a significant improvement in efficiency is achieved \cite{Panchi2022DCDS-B}. However, all of these schemes only have the first-order accuracy in time.
As shown in \Cref{fig:macromagnetics-simulations}, we depict the states of the system at $T = 2\; \mathrm{ns}$ and record the energy behaviors for different damping values, in which we fix the time step size $\Delta t = 1\; \mathrm{ps}$. The results demonstrate that our approach overcomes the dependence on the damping parameter as reported in \cite{garcia2003improved}. These demonstrate that our approach is unconditionally stable with respect to the damping parameter and reaches the optimal efficiency in micromagnetics simulations.

\section{Conclusions}
In this paper, we develop two Gauss-Seidel projection methods (GSPMs) with second-order accuracy for the Landau-Lifshitz (LL) equation. In the first-order GSPMs, heat equations are solved to construct the first-order approximations for the Laplacian operator in fractional steps. Instead, the biharmonic equation is solved to construct the second-order approximation for the Laplacian operator. Consequently, a combination of the semi-implicit second-order backward difference formula and the biharmonic equation yields the second-order consistent scheme for the LL equation. By applying the Gauss-Seidel manner and the projection step, we obtain the unconditional stable Scheme A, and conditionally stable Scheme B with a Courant-Friedrichs-Lewy constant of 0.25. We conduct numerical experiments to validate the approximation error. In addition, we present an unconditionally stable method for micromagnetics simulations. Numerical results showcase that our approach offers an accurate and efficient solution for simulating LL dynamics.

\section*{Acknowledgments}P. Li thanks Prof. Jingrun Chen for the fruitful discussions. X.-P. Wang acknowledges support from the National Natural Science Foundation of China (NSFC) (No. 12271461), the key project of NSFC (No. 12131010), the Hetao Shenzhen-Hong Kong Science and Technology Innovation Cooperation Zone Project (No.HZQSWSKCCYB-
2024016)


\bibliographystyle{amsplain}
\bibliography{refs}

\providecommand{\bysame}{\leavevmode\hbox to3em{\hrulefill}\thinspace}
\providecommand{\MR}{\relax\ifhmode\unskip\space\fi MR }
\providecommand{\MRhref}[2]{%
  \href{http://www.ams.org/mathscinet-getitem?mr=#1}{#2}
}
\providecommand{\href}[2]{#2}
\begin{thebibliography}{10}

\bibitem{Abert2019}
C.~Abert, \emph{Micromagnetics and spintronics: models and numerical methods},
  Eur. Phys. J. B \textbf{92} (2019), no.~6, 120.

\bibitem{akrivis2021higher}
G.~Akrivis, M.~Feischl, B.~Kov{\'a}cs, and C.~Lubich, \emph{Higher-order
  linearly implicit full discretization of the {Landau-Lifshitz-Gilbert}
  equation}, Math. Comput. \textbf{90} (2021), no.~329, 995--1038.

\bibitem{ALOUGES20121345}
F.~Alouges, E.~Kritsikis, and J.-C. Toussaint, \emph{A convergent finite
  element approximation for {Landau-Lifschitz-Gilbert} equation}, Physica B:
  Condens. Matter \textbf{407} (2012), no.~9, 1345--1349, 8th International
  Symposium on Hysteresis Modeling and Micromagnetics (HMM 2011).

\bibitem{doi:10.1137/20M1335431}
R.~An, H.~Gao, and W.~Sun, \emph{Optimal error analysis of euler and
  {Crank-Nicolson} projection finite difference schemes for {Landau-Lifshitz}
  equation}, SIAM J. Numer. Anal. \textbf{59} (2021), no.~3, 1639--1662.

\bibitem{10.1007/978-3-540-31852-1_17}
L.~Ba{\v{n}}as, \emph{Numerical methods for the {Landau-Lifshitz-Gilbert}
  equation}, Numerical Analysis and Its Applications (Berlin, Heidelberg)
  (Zhilin Li, Lubin Vulkov, and Jerzy Wa{\'{s}}niewski, eds.), Springer Berlin
  Heidelberg, 2005, pp.~158--165.

\bibitem{doi:10.1137/050631070}
S.~Bartels and A.~Prohl, \emph{Convergence of an implicit finite element method
  for the {Landau-Lifshitz-Gilbert} equation}, SIAM J. Numer. Anal. \textbf{44}
  (2006), no.~4, 1405--1419.

\bibitem{10.1145/79505.79507}
J.~R. Cash and Alan~H. Karp, \emph{A variable order {Runge-Kutta} method for
  initial value problems with rapidly varying right-hand sides}, ACM Trans.
  Math. Softw. \textbf{16} (1990), no.~3, 201–222.

\bibitem{doi:10.1137/22M1501143}
Q.~Cheng and J.~Shen, \emph{Length preserving numerical schemes for
  {Landau-Lifshitz} equation based on {Lagrange} multiplier approaches}, SIAM
  J. Sci. Comput. \textbf{45} (2023), no.~2, A530--A553.

\bibitem{10.1093/imanum/dri011}
I.~Cimr\'ak, \emph{Error estimates for a semi-implicit numerical scheme solving
  the {Landau-Lifshitz} equation with an exchange field}, IMA J. Numer. Anal.
  \textbf{25} (2005), no.~3, 611--634.

\bibitem{Cimrak2007}
I.~{}Cimr\'ak, \emph{A survey on the numerics and computations for the
  {Landau-Lifshitz} equation of micromagnetism}, Arch. Comput. Methods Eng.
  \textbf{15} (2008), 277--309.

\bibitem{Fuwa2012}
A.~Fuwa, T.~Ishiwata, and M.~Tsutsumi, \emph{Finite difference scheme for the
  {Landau-Lifshitz} equation}, Japan J. Indust. Appl. Math. \textbf{29} (2012),
  no.~1, 83--110.

\bibitem{garcia2003improved}
C.~J. Garc{\'\i}a-Cervera and W.~E, \emph{Improved {Gauss-Seidel} projection
  method for micromagnetics simulations}, IEEE Trans. Magn. \textbf{39} (2003),
  no.~3, 1766--1770.

\bibitem{Gilbert1955}
T.~Gilbert, \emph{A {Lagrangian} formulation of gyromagnetic equation of the
  magnetization field}, Phys. Rev. \textbf{100} (1955), 1243--1255.

\bibitem{gui2024thirdorder}
Y.~Gui, R.~Du, and C.~Wang, \emph{A third-order implicit-explicit {Runge-Kutta}
  method for {Landau-Lifshitz} equation with arbitrary damping parameters},
  2024.

\bibitem{SIAMRev2006DevelopmentLLG}
M.~Kru\v{z}\'ik and A.~Prohl, \emph{Recent developments in the modeling,
  analysis, and numerics of ferromagnetism}, SIAM Rev. \textbf{48} (2006),
  no.~3, 439--483.

\bibitem{LandauLifshitz1935}
L.~Landau and E.~Lifshitz, \emph{On the theory of the dispersion of magetic
  permeability in ferromagnetic bodies}, Phys. Z. Sowjetunion \textbf{8}
  (1935), 153--169.

\bibitem{li2023micromagnetics}
P.~Li, S.~Gu, J.~Lan, J.~Chen, W.~Ren, and R.~Du, \emph{Micromagnetics
  simulations and phase transitions of ferromagnetics with
  {Dzyaloshinskii-Moriya} interaction}, Commun. Nonlinear. Sci. Numer. Simul.
  \textbf{126} (2023), 107512.

\bibitem{Panchi2022DCDS-B}
P.~Li, Z.~Ma, R.~Du, and J.~Chen, \emph{A {Gauss-Seidel} projection method with
  the minimal number of updates for the stray field in micromagnetics
  simulations}, Discrete Contin. Dyn. Syst. Ser. B \textbf{27} (2022), no.~11,
  6401--6416.

\bibitem{LI2020109046}
P.~Li, C.~Xie, R.~Du, J.~Chen, and X.-P. Wang, \emph{Two improved
  {Gauss-Seidel} projection methods for {Landau-Lifshitz-Gilbert} equation}, J.
  Comput. Phys. \textbf{401} (2020), 109046.

\bibitem{ROMEO2008464}
A.~Romeo, G.~Finocchio, M.~Carpentieri, L.~Torres, G.~Consolo, and
  B.~Azzerboni, \emph{A numerical solution of the magnetization reversal
  modeling in a permalloy thin film using fifth order {Runge-Kutta} method with
  adaptive step size control}, Physica B: Condens. Matter \textbf{403} (2008),
  no.~2, 464--468.

\bibitem{WANG2001357}
X.-P. Wang, C.~J. Garc\'ia-Cervera, and W.~E, \emph{A {Gauss-Seidel} projection
  method for micromagnetics simulations}, J. Comput. Phys. \textbf{171} (2001),
  no.~1, 357--372.

\bibitem{xie2020second}
C.~Xie, C.~J. Garc{\'\i}a-Cervera, C.~Wang, Z.~Zhou, and J.~Chen,
  \emph{Second-order semi-implicit projection methods for micromagnetics
  simulations}, J. Comput. Phys. \textbf{404} (2020), 109104.

\end{thebibliography}
\end{document}